\documentclass[11pt,reqno]{amsart}
\usepackage{graphicx}
\usepackage{verbatim}
\usepackage{textcomp}
\usepackage{amssymb}
\usepackage{cite}
\usepackage{amsmath}
\usepackage{latexsym}
\usepackage{amscd}
\usepackage{amsthm}
\usepackage{mathrsfs}
\usepackage{xypic}
\usepackage{bm}
\usepackage{url}
\usepackage{hyperref}

\newtheorem{thm}{Theorem}[section]
\newtheorem{corr}[thm]{Corollary}
\newtheorem{lem}[thm]{Lemma}
\newtheorem{prop}[thm]{Proposition}

\theoremstyle{definition}

\theoremstyle{remark}
\newtheorem{rem}{Remark}[section]
\numberwithin{equation}{section}
\begin{document}
\title[Monotonicity formulas of eigenvalues and energy functionals
] {Monotonicity formulas of eigenvalues and energy functionals along
the rescaled List's extended Ricci flow}
\author{Guangyue Huang}
\address{College of Mathematics and Information Science, Henan Normal
University, Xinxiang, Henan 453007, People's Republic of China}
\email{hgy@henannu.edu.cn }
\author{Zhi Li}
\address{College of Mathematics and Information Science, Henan Normal
University, Xinxiang, Henan 453007, People's Republic of China}
\email{lizhihnsd@126.com} 
\thanks{The research of the first author was supported by NSFC No. 11371018, 11171091. }
\maketitle

\maketitle

\begin{abstract}
In this paper, we study monotonicity formulas of eigenvalues and
entropies along the rescaled List's extended Ricci flow. We derive
some monotonicity formulas of eigenvalues of Laplacian which
generalize those of Li in \cite{Liyi2014} and Cao-Hou-Ling in
\cite{Caohou2012}. Moreover, we also consider monotonicity formulas
of $\mathcal{F}_k$-functional which can be seen as a generalized
$\mathcal{F}$-functional corresponding with steady Ricci breathers,
and $\mathcal{W}_k$-functional which generalizes
$\mathcal{W}$-functional corresponding with expanding Ricci
breathers.
\end{abstract}

\section{Introduction}

Let $(M^n, g(t))$ be a compact Riemannian manifold, $g(t)$  be a
solution to the following List's extended Ricci flow which was
introduced by B. List:
\begin{equation}\label{Int1}\left\{\begin{array}{ll}
\frac{\partial}{\partial t}g =-2{\rm Ric}+2\alpha d\varphi\otimes d\varphi,\\
\varphi_t=\Delta \varphi,
\end{array}\right.
\end{equation} where $\alpha>0$ is a real constant, $\varphi=\varphi(t)$ is
a smooth scalar function defined on $M^n$ and $\Delta$ denotes the
Laplacian given by $g(t)$. When $\alpha=2$ and
$\alpha=\frac{n-1}{n-2}$, the extended Ricci flow \eqref{Int1} have
been studied by List in \cite{List2006} and \cite{List2008},
respectively. Denote by $S_{ij}=R_{ij}-\alpha\varphi_i\varphi_j$ a
symmetric two-tensor. Then \eqref{Int1} becomes
\begin{equation}\label{Int2}\left\{\begin{array}{ll}
\frac{\partial}{\partial t}g_{ij}=-2S_{ij},\\
\varphi_t=\Delta \varphi.
\end{array}\right.
\end{equation}

In this paper, we consider the following rescaled List's extended
Ricci flow
\begin{equation}\label{Int3}\left\{\begin{array}{ll}
\frac{\partial}{\partial t}g_{ij}=-2(S_{ij}-\frac{r}{n}g_{ij}),\\
\varphi_t=\Delta \varphi,
\end{array}\right.
\end{equation} where $r=r(t)$ is a function depending only on $t$.
In particular, \eqref{Int2} can be seen as a special case of
\eqref{Int3} when $r=0$. On the other hand, if
$r(t)=(\int_MS\,dv)/(\int_M\,dv)$ (that is, $r$ is the average value
of $S$), then \eqref{Int3} can be seen as the extended Hamilton
normalized flow under the List's extended Ricci flow.  Here
\begin{equation}\label{Int4}S=g^{ij}S_{ij}=R-\alpha|\nabla\varphi|^2\end{equation}
is the trace of the two-tensor $S_{ij}$. For a constant $k\geq1$, we
defined the $\mathcal{F}_k$-functional as follows
\begin{equation}\label{Int5}\mathcal{F}_k=\int\limits_M(|\nabla
f|^2+kS)\,e^{-f}dv.
\end{equation}

If we define
\begin{equation}\label{eigSec1}
\lambda(g)=\inf\limits_f\mathcal{F}_k,
\end{equation} where the infimum is taken over all smooth function $f$ which satisfies
\begin{equation}\label{eigSec2}\int\limits_Me^{-f}\,dv=1,\end{equation} then
the nondecreasing of the $\mathcal{F}_k$-functional implies the
nondecreasing of $\lambda(g)$. In particular, $\lambda(g)$ defined
in \eqref{eigSec1} is the lowest eigenvalue of the operator
$-4\Delta+kS$. In the first part, we consider eigenvalues of the
operator
\begin{equation}\label{eigSec3}-\Delta+bS\end{equation} with $b$ a constant. For $b=0$, we
first derive the following evolution equation of eigenvalues on
Laplacian under the rescaled List's extended Ricci flow
\eqref{Int3}. That is, we obtain

\begin{thm}\label{4thm1}
Let $\lambda^{-\Delta}(t)$ be the eigenvalue of the operator
$-\Delta$ corresponding to the normalized eigenfunction $u$, that
is,
$$-\Delta
u=\lambda^{-\Delta} u,\ \ \ \  \int\limits_Mu^2\,dv=1.$$ Then under
the rescaled List's extended Ricci flow \eqref{Int3},
\begin{equation}\label{4Eig10}\aligned
\frac{d}{d t}\lambda^{-\Delta}=&-\frac{2r}{n}\lambda^{-\Delta}
+\int\limits_M\Big(\lambda^{-\Delta}Su^2-S|\nabla u|^2
+2S^{ij}u_{i}u_j\Big)\,dv.
\endaligned\end{equation}
\end{thm}


\begin{thm}\label{4thm2}
Let $\lambda^{-\Delta+\frac{1}{2}S}(t)$ be the eigenvalue of the
operator $-\Delta+\frac{1}{2}S$ corresponding to the normalized
eigenfunction $u$, that is,
$$(-\Delta+\frac{1}{2}S)
u=\lambda^{-\Delta+\frac{1}{2}S} u,\ \ \ \ \int\limits_Mu^2\,dv=1.$$
Then under the rescaled List's extended Ricci flow \eqref{Int3},
\begin{equation}\label{4Eig11}\aligned
\frac{d}{d
t}\lambda^{-\Delta+\frac{1}{2}S}=&-\frac{2r}{n}\lambda^{-\Delta+\frac{1}{2}S}
+\int\limits_M\Big[|S_{ij}|^2u^2+2S^{ij}u_{i}u_j+\alpha(\Delta\varphi)^2u^2\Big]\,dv.
\endaligned\end{equation}
Moreover, if $S_{ij}(t)\geq0$ for all $t$, the eigenvalues  of the
operator $-\Delta+\frac{1}{2}S$ satisfy
\begin{equation}\label{4Eig13}\aligned
\frac{d}{d
t}\Big(\lambda^{-\Delta+\frac{1}{2}S}e^{\frac{2}{n}\int_0^tr(s)\,ds}\Big)=&
e^{\frac{2}{n}\int_0^tr(s)\,ds}\Bigg\{\int\limits_M\Big[|S_{ij}|^2u^2+2S^{ij}u_{i}u_j\\
&+\alpha(\Delta\varphi)^2u^2\Big]\,dv\Bigg\} \geq0
\endaligned\end{equation} and
$\lambda^{-\Delta+\frac{1}{2}S}e^{\frac{2}{n}\int_0^tr(s)\,ds}$ is
nondecreasing under the rescaled List's extended Ricci flow
\eqref{Int3}. Furthermore, the monotonicity is strict unless the
metric is Ricci flat.

\end{thm}

It is well-known that, under the List's extended Ricci flow
\eqref{Int2},  the nonnegativity of $S$ is preserved. In this paper,
we will prove that the nonnegativity of $S$ is also preserved under
the rescaled List's extended Ricci flow \eqref{Int3} for all $r(t)$.
That is,

\begin{thm}\label{nonnegativity-thm1}
The nonnegativity of $S$ is preserved under the rescaled List's
extended Ricci flow \eqref{Int3}.

\end{thm}

In order to state the following results on eigenvalues, we first
introduce the following definition: $$S_{\min}(0)=\min\limits_{x\in
M} S(x,0).$$ By virtue of Theorem \ref{nonnegativity-thm1}, we prove
the following

\begin{thm}\label{spec4thm1}
Let $(g(t),\varphi(t))$ be a solution to the rescaled List's
extended Ricci flow \eqref{Int3} with $S_{ij}(t)\geq \theta
g_{ij}(t)$ holding for some $\theta\geq\frac{1}{2}$. Let
$\lambda^{-\Delta}(t)$ be the eigenvalue of the operator $-\Delta$.

(1) If $S_{\min}(0)\geq0$, then
$\lambda^{-\Delta}e^{\frac{2}{n}\int_0^tr(s)\,ds}$ is nondecreasing
along the rescaled List's extended Ricci flow \eqref{Int3}.

(2) For all $t$, we have
\begin{equation}\label{spec4thm1-1} \frac{d}{d
t}\ln\Big(\lambda^{-\Delta}e^{\frac{2}{n}\int_0^tr(s)\,ds}\Big)\geq
2\theta x(t).
\end{equation} Moreover,
$\lambda(t)$ has the lower bound
\begin{equation}\label{spec4thm1-2} \lambda^{-\Delta}(t)e^{\frac{2}{n}\int_0^tr(s)\,ds}\geq
\lambda(0)e^{2\theta \int_0^tx(s)\,ds}
\end{equation} depending only on $t$, where
$$x(t)=\frac{S_{\min}(0)e^{-\frac{2}{n}\int_0^t
r(s)ds}}{1-\frac{2}{n}S_{\min}(0)\int_0^t\Big(e^{-\frac{2}{n}\int_0^s
r(\tilde{s})d\tilde{s}}\Big)ds}.$$

\end{thm}

In particular, for compact Riemannian surfaces, we obtain the
following consequences from Theorem \ref{spec4thm1}:

\begin{corr}\label{spec4thm2}
Let $(g(t),\varphi(t))$ be a solution to the rescaled List's
extended Ricci flow \eqref{Int3} on $M^2$. Let
$\lambda^{-\Delta}(t)$ be the eigenvalue of the operator $-\Delta$.

(i) If $$R_{ij}\leq\epsilon u_iu_j,$$ where
$\epsilon\leq\frac{2\alpha(\theta-1)}{2\theta-1}$ with
$\theta>\frac{1}{2}$, then the following holds:

(1) If $S_{\min}(0)\geq0$, then
$\lambda^{-\Delta}e^{\int_0^tr(s)\,ds}$ is nondecreasing along the
rescaled List's extended Ricci flow \eqref{Int3}.

(2) For all $t$, we have
\begin{equation}\label{spec4thm1-3} \frac{d}{d
t}\ln\Big(\lambda^{-\Delta}e^{\int_0^tr(s)\,ds}\Big)\geq 2\theta
x(t).
\end{equation} Moreover,
$\lambda(t)$ has the lower bound
\begin{equation}\label{spec4thm1-4} \lambda^{-\Delta}(t)e^{\int_0^tr(s)\,ds}\geq
\lambda(0)e^{2\theta\int_0^t x(s)\,ds},
\end{equation} depending only on $t$, where
$$x(t)=\frac{S_{\min}(0)e^{-\int_0^t
r(s)ds}}{1-S_{\min}(0)\int_0^t\Big(e^{-\int_0^s
r(\tilde{s})d\tilde{s}}\Big)ds}.$$

(ii) If $|\nabla \varphi|^2g_{ij}\geq2\varphi_i\varphi_j$, then the
following holds:

(1) If $S_{\min}(0)\geq0$, then
$\lambda^{-\Delta}e^{\int_0^tr(s)\,ds}$ is nondecreasing along the
rescaled List's extended Ricci flow \eqref{Int3}.

(2) For all $t$, we have
\begin{equation}\label{spec4thm1-3} \frac{d}{d
t}\ln\Big(\lambda^{-\Delta}e^{\int_0^tr(s)\,ds}\Big)\geq x(t).
\end{equation} Moreover,
$\lambda(t)$ has the lower bound
\begin{equation}\label{spec4thm1-4} \lambda^{-\Delta}(t)e^{\int_0^tr(s)\,ds}\geq
\lambda(0)e^{\int_0^tx(s)\,ds}
\end{equation} depending only on $t$, where
$$x(t)=\frac{S_{\min}(0)e^{-\int_0^t
r(s)ds}}{1-S_{\min}(0)\int_0^t\Big(e^{-\int_0^s
r(\tilde{s})d\tilde{s}}\Big)ds}.$$
\end{corr}

\begin{rem}
It should be pointed out that for $r=0$ and $\alpha=2$, our above
results on eigenvalues reduce to the corresponding results of Li in
\cite{Liyi2014}. In particular, our Theorem \ref{nonnegativity-thm1}
is new.

\end{rem}

\begin{rem}
Some related results for monotonicity formulas of eigenvalues on
Laplacian along the Ricci flow, we refer to
\cite{Ma2006,Cao2007,Cao2008} and among others \cite{Fang2015} for
later development.
\end{rem}

Next, we study monotonicity formulas of eigenvalues on Laplacian on
Riemannian surfaces. We obtain the following results:

\begin{thm}\label{twodim4thm3}

Let $\lambda^{-\Delta+bS}(t)$ be the eigenvalue of the operator
$-\Delta+bS$ with normalized eigenfunction $u$ on $M^2$ with
$S(t)\geq0$ holding for all $t$, that is,
$$(-\Delta+bS)
u=\lambda^{-\Delta+bS} u,\ \ \ \ \int\limits_Mu^2\,dv=1.$$ Then if
$|\nabla \varphi|^2g_{ij}\geq2\varphi_i\varphi_j$ and
$0<b\leq\frac{1}{2}$, we have
\begin{equation}\label{4Eig11}\aligned
\frac{d}{d t}\Big(\lambda^{-\Delta+bS}e^{\int_0^tr(s)\,ds}\Big)\geq&
e^{\int_0^tr(s)\,ds}\Bigg\{\int\limits_M\Big\{2b^2S^2u^2+(1-2b)\lambda
Su^2+2bS|\nabla u|^2\\
&+b\alpha^2|\nabla \varphi|^4u^2+2b\alpha(\Delta\varphi)^2u^2]\Big\}\,dv\Bigg\}\\
&\geq0.
\endaligned\end{equation}
\end{thm}

We also obtain the following bounds for eigenvalues of the operator
$-\Delta+bS$ on compact Riemannian surfaces.

\begin{thm}\label{two4thm4}

Let $\lambda^{-\Delta+bS}(t)$ be the eigenvalue of the operator
$-\Delta+bS$  on $M^2$. If $|\nabla
\varphi|^2g_{ij}\geq2\varphi_i\varphi_j$ and $0<b\leq\frac{1}{2}$,
we have
\begin{equation}\label{two-eigenvalue9} [1-tS_{\min}(0)]\lambda
-\frac{b^2S^2_{\min}(0)}{2}\ln[1-tS_{\min}(0)]
\end{equation} is nondecreasing under
the List's extended Ricci flow \eqref{Int2}. Moreover, $\lambda(t)$
has the lower bound
\begin{equation}\label{two-eigenvalue10} \lambda(t)\geq\frac{1}{1-tS_{\min}(0)}\lambda(0)
+\frac{b^2S_{\min}(0)}{2[1-tS_{\min}(0)]}\ln[1-tS_{\min}(0)]
\end{equation} depending only on $t$.

\end{thm}

\begin{thm}\label{two4thm5}

Let $\lambda^{-\Delta+bS}(t)$ be the eigenvalue of the operator
$-\Delta+bS$  on $M^2$ with $r>0$ and $S_{\min}(0)>0$. If $|\nabla
\varphi|^2g_{ij}\geq2\varphi_i\varphi_j$ and $0<b\leq\frac{1}{2}$,
we have
\begin{equation}\label{two-eigenvalue17} \frac{d}{d
t}(\ln\lambda)\geq \frac{S_{\min}(0)e^{-\int_0^t
r(s)ds}}{1-S_{\min}(0)\int_0^t\Big(e^{-\int_0^s
r(\tilde{s})d\tilde{s}}\Big)ds}-r
\end{equation} under
the rescaled List's extended Ricci flow \eqref{Int3}. Moreover,
$\lambda(t)$ has the lower bound
\begin{equation}\label{two-eigenvalue10} \lambda(t)\geq
\lambda(0)e^{\int_0^t\tilde{x}(s)\,ds},
\end{equation} depending only on $t$, where $$\tilde{x}(t)=\frac{S_{\min}(0)e^{-\int_0^t
r(s)ds}}{1-S_{\min}(0)\int_0^t\Big(e^{-\int_0^s
r(\tilde{s})d\tilde{s}}\Big)ds}-r.$$

\end{thm}

\begin{rem}
When $r=0$, Theorem \ref{twodim4thm3} becomes Theorem 1.6  of
Cao-Hou-Ling in \cite{Caohou2012}. When $\varphi=0$ and
$r(t)=(\int_MR\,dv)/(\int_M\,dv)$, our Theorems \ref{two4thm4},
\ref{two4thm5} reduce to Theorem 3.4 and Theorem 3.3 of Cao-Hou-Ling
in \cite{Caohou2012}, respectively.

\end{rem}

In the rest of this paper, we consider monotonicity formulas of
$\mathcal{F}_k$-functional which can be seen as a generalized
$\mathcal{F}$-functional corresponding with steady Ricci breathers,
and $\mathcal{W}_k$-functional which can be seen as a generalized
$\mathcal{W}$-functional corresponding with expanding Ricci
breathers. Under the following evolution equation
\begin{equation}\label{2Proof1}\left\{\begin{array}{ll}
\frac{\partial}{\partial t}g_{ij}=-2(S_{ij}-\frac{r}{n}g_{ij}),\\
\varphi_t=\Delta \varphi,\\
f_t=-\Delta f+|\nabla f|^2-S+r,
\end{array}\right.\end{equation}
we proved the following results:

\begin{thm}\label{2thm1} Under the system \eqref{2Proof1}, we have
\begin{equation}\label{add2Proof17}\aligned
\frac{d}{dt}\mathcal{F}_k=&-\frac{2r}{n}\mathcal{F}_k+2(k-1)\int\limits_M\Big(|S_{ij}|^2+\alpha(\Delta
\varphi)^2\Big)\,e^{-f}dv\\
&+2\int\limits_M\Big(|S_{ij}+f_{ij}|^2+\alpha|\Delta
\varphi-\langle\nabla f,\nabla\varphi\rangle|^2\Big)\,e^{-f}dv,
\endaligned\end{equation}
or equivalently,
\begin{equation}\label{add2Proof18}\aligned
\frac{d}{dt}\mathcal{F}_k=&\frac{2r}{n}(\mathcal{F}_k-kr)
+2(k-1)\int\limits_M\Big(|S_{ij}-\frac{r}{n}g_{ij}|^2+\alpha(\Delta
\varphi)^2\Big)\,e^{-f}dv\\
&+2\int\limits_M\Big(|S_{ij}+f_{ij}-\frac{r}{n}g_{ij}|^2+\alpha|\Delta
\varphi-\langle\nabla f,\nabla\varphi\rangle|^2\Big)\,e^{-f}dv.
\endaligned\end{equation}

\end{thm}

\begin{rem}

It was pointed out that for $\alpha=2$ and $k=1$, List
\cite{List2006} and M\"{u}ller \cite{Muller2012} studied the
monotonicity of $\mathcal{F}_k$-functional under the List's extended
Ricci flow \eqref{Int2}. Later, In \cite{Liyi2014}, Li studied the
monotonicity of $\mathcal{F}_k$-functional for $\alpha=2$ and all
$k\geq1$ under \eqref{Int2}. In particular, when $\varphi=0$ and
$k=1$, the $\mathcal{F}_k$-functional \eqref{Int5} becomes the
Perelman's $\mathcal{F}$-functional.

\end{rem}

Applying \eqref{add2Proof18} in Theorem \ref{2thm1}, we can also
obtain the following result:

\begin{thm}\label{1th1}

Let $(M^n,g(t))$ be a compact Riemannian manifold with $g(t)$
satisfying the rescaled List's extended Ricci flow \eqref{Int3}. We
let $\lambda(t)$ be the lowest eigenvalue of the operator
$-4\Delta+kS$ with $k\geq1$. If the average value of $S$ is
nonnegative for all $t$, then $\lambda(t)$ is nondecreasing along
\eqref{Int3}. Moreover, the monotonicity is strict unless the metric
is Einstein.

\end{thm}

As a direct application of \eqref{add2Proof17} in Theorem
\ref{2thm1}, we can obtain the following results:

\begin{corr}\label{2corr1} Under the system \eqref{2Proof1} and $k\geq1$, we have
\begin{equation}\label{add2Proof19}\aligned
\frac{d}{d
t}\Big(\mathcal{F}_ke^{\frac{2}{n}\int_0^tr(s)\,ds}\Big)=&e^{\frac{2}{n}\int_0^tr(s)\,ds}\Bigg\{2(k-1)\int\limits_M\Big(|S_{ij}|^2+\alpha(\Delta
\varphi)^2\Big)\,e^{-f}dv\\
&+2\int\limits_M\Big(|S_{ij}+f_{ij}|^2+\alpha|\Delta
\varphi-\langle\nabla
f,\nabla\varphi\rangle|^2\Big)\,e^{-f}dv\Bigg\}\\
\geq&0,
\endaligned\end{equation} which shows that
$\mathcal{F}_ke^{\frac{2}{n}\int_0^tr(s)\,ds}$ is nondecreasing
along \eqref{2Proof1}. Moreover, the monotonicity is strict unless
the metric is Ricci flat.

\end{corr}

\begin{rem} Choosing $r=0$ and $\alpha=2$, then \eqref{add2Proof17} reduces the formula (1-6)
in Theorem 1.1 of Li \cite{Liyi2014}. On the other hand, under the
normalized Ricci flow of Hamilton, Li in \cite{Lijunfang2011} also
obtained a similar result as Theorem \ref{1th1}.
\end{rem}

As in \cite{Liyi2014}, for a constant $k\geq1$, we define the
following $\mathcal{W}_k$-functional:
\begin{equation}\label{Int10}
\mathcal{W}_k=\tau^2\int\limits_M[k(S+\frac{n}{2\tau})+|\nabla
f|^2]e^{-f}\,dv.
\end{equation}
Under the following coupled system
\begin{equation}\label{3Ent1}\left\{\begin{array}{ll}
\frac{\partial}{\partial t}g_{ij}=-2(S_{ij}-\frac{r}{n}g_{ij}),\\
\varphi_t=\Delta \varphi,\\
f_t=-\Delta f+|\nabla f|^2-S+r\\
\tau_t=1,
\end{array}\right.
\end{equation} we obtain the following results:

\begin{thm}\label{3thm1} Under the system \eqref{3Ent1} and $k\geq1$, we have

\begin{equation}\label{3Ent2}\aligned
\frac{d}{d
t}\mathcal{W}_k=&2\tau^2\Bigg\{-\frac{r}{n}\mathcal{F}_k+(k-1)\int\limits_M\Big(|S_{ij}+\frac{1}{2\tau}g_{ij}|^2
+\alpha(\Delta \varphi)^2\Big)\,e^{-f}dv\\
&+\int\limits_M\Big(|S_{ij}+f_{ij}+\frac{1}{2\tau}g_{ij}|^2
+\alpha|\Delta \varphi-\langle\nabla
f,\nabla\varphi\rangle|^2\Big)\,e^{-f}dv\Bigg\},
\endaligned\end{equation}
or equivalently,
\begin{equation}\label{add3Ent2}\aligned
\frac{d}{d
t}\mathcal{W}_k=&2\tau^2\Bigg\{\frac{r}{n}(\mathcal{F}_k-kr)+\frac{kr}{\tau}\\
&+(k-1)\int\limits_M\Big(|S_{ij}+\frac{1}{2\tau}g_{ij}-\frac{r}{n}g_{ij}|^2
+\alpha(\Delta \varphi)^2\Big)\,e^{-f}dv\\
&+\int\limits_M\Big(|S_{ij}+f_{ij}+\frac{1}{2\tau}g_{ij}-\frac{r}{n}g_{ij}|^2
+\alpha|\Delta \varphi-\langle\nabla
f,\nabla\varphi\rangle|^2\Big)\,e^{-f}dv\Bigg\}.
\endaligned\end{equation}

\end{thm}

As a direct application of the formula \eqref{3Ent2}, we obtain

\begin{corr}\label{corr2} Under the system \eqref{3Ent1} and $k\geq1$, we have

\begin{equation}\label{corr1}\aligned
\frac{d}{d
t}\mathcal{W}_k+\frac{2r}{n}\tau^2\mathcal{F}_k=&2\tau^2\Bigg\{(k-1)\int\limits_M\Big(|S_{ij}+\frac{1}{2\tau}g_{ij}|^2
+\alpha(\Delta \varphi)^2\Big)\,e^{-f}dv\\
&+\int\limits_M\Big(|S_{ij}+f_{ij}+\frac{1}{2\tau}g_{ij}|^2
+\alpha|\Delta \varphi-\langle\nabla
f,\nabla\varphi\rangle|^2\Big)\,e^{-f}dv\Bigg\}\\
\geq&0.
\endaligned\end{equation} Moreover, the equality in \eqref{corr1}
holds if and only if the metric is Einstein.

\end{corr}

\begin{rem} When $r=0$ and $\alpha=2$, then \eqref{3Ent2} reduces the formula
(1-8) in Theorem 1.3 of Li \cite{Liyi2014}.
\end{rem}


This paper is organized as follows: In Section two, we study
monotonicity formulas for eigenvalues of the operator $-\Delta+bS$.
We mainly consider monotonicity formulas for eigenvalues of the
Laplacian under the rescaled List's extended Ricci flow \eqref{Int3}
and give some dependent lower bounds of eigenvalues. Moreover,
Theorem \ref{4thm1}-Theorem \ref{spec4thm1} have been proved in this
part. In Section three, we deal with monotonicity formulas for
eigenvalues of the operator $-\Delta+bS$ with $0<b\leq\frac{1}{2}$
on compact Riemannian surfaces and obtain some interesting results.
In this part, we give proofs of Theorem \ref{twodim4thm3}-Theorem
\ref{two4thm5}. In the last Section four, we study the first
variations of $\mathcal{F}_k$-functional and
$\mathcal{W}_k$-functional. These functionals are very important to
study entropies corresponding to the List's extended Ricci flow
\eqref{Int2}. In this part, Theorem \ref{2thm1}-Theorem \ref{3thm1}
have been proved.

\section{Proof of Theorem \ref{4thm1}-Theorem \ref{spec4thm1} }

We first give a Lemma which will be used later.

\begin{lem}\label{2eigenvaluelem1} Under the rescaled List's extended Ricci flow
\eqref{Int3}, we have
\begin{equation}\label{2eigenvalue1}
\frac{\partial}{\partial
t}g^{ij}=2g^{ik}g^{jl}(S_{kl}-\frac{r}{n}g_{kl}),
\end{equation}
\begin{equation}\label{2eigenvalue2}
(dv)_t=(-S+r)\,dv,
\end{equation}
\begin{equation}\label{2eigenvalue3}
S_t=\Delta S+2|S_{ij}|^2-\frac{2r}{n}S+2\alpha(\Delta\varphi)^2.
\end{equation} In particular, when $r(t)=(\int_MS\,dv)/(\int_M\,dv)$ (that is, $r$ is the average value
of $S$), then \eqref{2eigenvalue2} shows that the volume of
$(M,g(t))$ is a constant for all $t$.

\end{lem}

\proof It has been shown in \cite{List2006}(see Lemma 1.4 in
\cite{List2006}) that if the metric $g(t)$ satisfies
\begin{equation}\label{2eigenvalue4}
\frac{\partial}{\partial t}g_{ij}=v_{ij},
\end{equation} where $v_{ij}$ is a symmetric two tensor, then
\begin{equation}\label{2eigenvalue5}
\frac{\partial}{\partial t}g^{ij}=-v^{ij},
\end{equation}
\begin{equation}\label{2eigenvalue6}
(dv)_t=\frac{1}{2}({\rm tr_g}v)\,dv
\end{equation} with ${\rm tr_g}v$ denoting the trace of $v_{ij}$
with respect to $g$, and
\begin{equation}\label{2eigenvalue7}
R_t=-\Delta({\rm
tr_g}v)+g^{kj}g^{il}v_{ij,kl}-g^{kj}g^{il}R_{kl}v_{ij}.
\end{equation}
Hence, inserting $v_{ij}=-2(S_{ij}-\frac{r}{n}g_{ij})$, ${\rm
tr_g}v=-2(S-r)$ into \eqref{2eigenvalue5} and \eqref{2eigenvalue6},
we obtain \eqref{2eigenvalue1} and \eqref{2eigenvalue2},
respectively.

Next, we prove \eqref{2eigenvalue3}. In fact, from
\eqref{2eigenvalue7}, we have
\begin{equation}\label{2eigenvalue8}\aligned
R_t=&2\Delta
S-2g^{kj}g^{il}S_{ij,kl}+2g^{kj}g^{il}R_{kl}(S_{ij}-\frac{r}{n} g_{ij})\\
=&2\Delta
S-2g^{kj}g^{il}S_{ij,kl}+2g^{kj}g^{il}(S_{kl}+\alpha \varphi_k\varphi_l)(S_{ij}-\frac{r}{n} g_{ij})\\
=&2\Delta S-2g^{kj}g^{il}S_{ij,kl}+2|S_{ij}|^2-\frac{2r}{n}S+2\alpha
S_{ij}\varphi^i\varphi^j-\alpha\frac{2r}{n}|\nabla \varphi|^2.
\endaligned\end{equation}
By the definition of $S_{ij}$ and the contracted Bianchi indentity,
we also have
\begin{equation}\label{2eigenvalue9}\aligned g^{kj}S_{ij,k}=&g^{kj}R_{ij,k}
-g^{kj}\alpha(\varphi_i\varphi_j)_{,k}\\
=&\frac{1}{2}(R-\alpha|\nabla\varphi|^2)_{,i}-\alpha(\Delta\varphi)\varphi_i\\
=&\frac{1}{2}S_{,i}-\alpha(\Delta\varphi)\varphi_i.
\endaligned\end{equation} Thus, \eqref{2eigenvalue8} becomes
\begin{equation}\label{2eigenvalue10}\aligned
R_t=&\Delta S+2|S_{ij}|^2-\frac{2r}{n}S+2\alpha
S_{ij}\varphi^i\varphi^j-\alpha\frac{2r}{n}|\nabla
\varphi|^2+2\alpha(\Delta\varphi)^2\\
&+2\alpha\langle\nabla\Delta\varphi,\nabla\varphi\rangle.
\endaligned\end{equation}
It follows that
\begin{equation}\label{2eigenvalue11}\aligned
S_t=&(R-\alpha|\nabla\varphi|^2)_t\\
=&\Delta S+2|S_{ij}|^2-\frac{2r}{n}S+2\alpha(\Delta\varphi)^2
\endaligned\end{equation}
and the desired \eqref{2eigenvalue3} follows. We complete the proof
of Lemma \ref{2eigenvaluelem1}.

Let $u$ be the eigenfunction corresponding to eigenvalue $\lambda$
of the operator $-\Delta+bS$, that is,
\begin{equation}\label{4Eig1}
(-\Delta+bS)u=\lambda u.
\end{equation} Multiplying both sides of \eqref{4Eig1} with $u$
and integrating on $M$, we have
\begin{equation}\label{4Eig2}
\lambda=\int\limits_M (|\nabla u|^2+bSu^2)\,dv.
\end{equation}
Using Lemma \ref{2eigenvaluelem1}, we have
\begin{equation}\label{4Eig3}\aligned
\frac{d}{d t}\lambda=&\int\limits_M
(2u_iu_j\frac{\partial }{\partial t}g^{ij}+2(u_t)^iu_i+bS_tu^2+2bSuu_t)\,dv\\
&+\int\limits_M (|\nabla u|^2+bSu^2)(-S+r)\,dv\\
=&\int\limits_M(2S^{ij}u_iu_j-\frac{2r}{n}|\nabla u|^2-2u_t\Delta u+bS_tu^2+2bSuu_t)\,dv\\
&+\int\limits_M (|\nabla u|^2+bSu^2)(-S+r)\,dv.
\endaligned\end{equation}
From \eqref{2eigenvalue9} and the Stokes formula, we have
\begin{equation}\label{4Eig4}\aligned
2\int\limits_MS^{ij}u_iu_j\,dv=\int\limits_M(-S_{,i}u^iu-2S^{ij}u_{ij}u
+2\alpha (\Delta\varphi)\langle\nabla u,\nabla\varphi\rangle u)\,dv.
\endaligned\end{equation} On the other
hand,
\begin{equation}\label{4Eig5}\aligned
-\int\limits_M|\nabla u|^2S\,dv=&\int\limits_M(S\Delta
u+S_{,i}u^i)u\,dv.
\endaligned\end{equation} Inserting \eqref{4Eig4} and \eqref{4Eig5}
into \eqref{4Eig3} yields
\begin{equation}\label{4Eig6}\aligned
\frac{d}{d t}\lambda=&\int\limits_M\Bigg\{-\frac{2r}{n}|\nabla
u|^2+bS_tu^2-2S^{ij}u_{ij}u+2\alpha(\Delta\varphi)\langle\nabla u,\nabla\varphi\rangle u\\
&+2u_t(-\Delta u+bSu)-Su(-\Delta u+bSu)+r(|\nabla u|^2+bSu^2)\Bigg\}\,dv\\
=&\int\limits_M\Bigg\{-\frac{2r}{n}|\nabla
u|^2+bS_tu^2-2S^{ij}u_{ij}u+2\alpha(\Delta\varphi)\langle\nabla u,\nabla\varphi\rangle u\\
&+\lambda2u_tu+\lambda(-S+r) u^2\Bigg\}\,dv\\
=&\int\limits_M\Bigg\{-\frac{2r}{n}|\nabla
u|^2+bS_tu^2-2S^{ij}u_{ij}u+2\alpha(\Delta\varphi)\langle\nabla u,\nabla\varphi\rangle u\Bigg\}\\
&+\lambda\Bigg(\int\limits_M u^2\,dv\Bigg)_t\\
=&\int\limits_M\Bigg\{-\frac{2r}{n}|\nabla
u|^2+bS_tu^2-2S^{ij}u_{ij}u+2\alpha(\Delta\varphi)\langle\nabla
u,\nabla\varphi\rangle u\Bigg\}\,dv.
\endaligned\end{equation}
Applying \eqref{2eigenvalue3} in Lemma \ref{2eigenvaluelem1} into
\eqref{4Eig6} yields
\begin{equation}\label{4Eig7}\aligned
\frac{d}{d t}\lambda=&\int\limits_M\Big\{-\frac{2r}{n}(|\nabla
u|^2+bSu^2)+2b|S_{ij}|^2u^2-2S^{ij}u_{ij}u+bu^2(\Delta S)\\
&+2b\alpha(\Delta\varphi)^2u^2+2\alpha(\Delta\varphi)\langle\nabla
u,\nabla\varphi\rangle u\Big\}\,dv\\
=&-\frac{2r}{n}\lambda+\int\limits_M\Big\{2b|S_{ij}|^2u^2-2S^{ij}u_{ij}u+bu^2\Delta S\\
&+2b\alpha(\Delta\varphi)^2u^2+2\alpha(\Delta\varphi)\langle\nabla
u,\nabla\varphi\rangle u\Big\}\,dv.
\endaligned\end{equation}
Using the Stokes formula again, we have
\begin{equation}\label{4Eig8}\aligned
-2\int\limits_MS^{ij}u_{ij}u\,dv=&\int\limits_M(2S^{ij}{}_{,j}u_{i}u\,dv
+2S^{ij}u_{i}u_j)\,dv\\
=&\int\limits_M[S_{,i}u_{i}u+2S^{ij}u_{i}u_j-2\alpha(\Delta\varphi)\langle\nabla
u,\nabla\varphi\rangle u]\,dv\\
=&\int\limits_M\Big[-\frac{1}{2}u^2\Delta
S+2S^{ij}u_{i}u_j-2\alpha(\Delta\varphi)\langle\nabla
u,\nabla\varphi\rangle u\Big]\,dv.
\endaligned\end{equation}
Therefore, \eqref{4Eig7} becomes
\begin{equation}\label{4Eig9}\aligned
\frac{d}{d t}\lambda=&-\frac{2r}{n}\lambda
+\int\limits_M\Big\{2b|S_{ij}|^2u^2+2S^{ij}u_{i}u_j+(b-\frac{1}{2})u^2\Delta
S+2b\alpha(\Delta\varphi)^2u^2\Big\}\,dv\\
=&-\frac{2r}{n}\lambda
+\int\limits_M\Big\{2b|S_{ij}|^2u^2+2S^{ij}u_{i}u_j+2b\alpha(\Delta\varphi)^2u^2\\
&+(2b-1)S[(bS-\lambda)u^2+|\nabla u|^2]\Big\}\,dv.
\endaligned\end{equation}

Hence, the following consequence follows:

\begin{prop}\label{Propth1}
Let $\lambda^{-\Delta+bS}(t)$ be the eigenvalue of the operator
$-\Delta+bS$  corresponding to the normalized eigenfunction $u$,
that is,
$$(-\Delta+bS)
u=\lambda^{-\Delta+bS} u,\ \ \ \  \int\limits_Mu^2\,dv=1.$$ Then
\begin{equation}\label{4Eig10}\aligned
\frac{d}{d t}\lambda^{-\Delta+bS}=&-\frac{2r}{n}\lambda^{-\Delta+bS}
+\int\limits_M\Big\{2b|S_{ij}|^2u^2+2S^{ij}u_{i}u_j+2b\alpha(\Delta\varphi)^2u^2\\
&+(2b-1)S[(bS-\lambda)u^2+|\nabla u|^2]\Big\}\,dv.
\endaligned\end{equation}

\end{prop}

\vspace*{2mm} \noindent{\bf Proof of Theorems \ref{4thm1} and
\ref{4thm2}.}  From Proposition \ref{Propth1}, it is easy to derive
Theorem \ref{4thm1} by letting $b=0$ and Theorem \ref{4thm2} by
letting $b=\frac{1}{2}$, respectively.

Next, we give the proof of Theorem \ref{nonnegativity-thm1} by the
Lemma 2.12 in \cite{Chow2007} (see Page 99 in \cite{Chow2007}).

\vspace*{2mm} \noindent{\bf Proof of Theorem
\ref{nonnegativity-thm1}.} Using the Cauchy inequality
$$|S_{ij}|^2\geq\frac{1}{n}S^2,$$ we obtain from \eqref{2eigenvalue3}
\begin{equation}\label{nonnegativity1}\aligned
S_t\geq& \Delta S+\frac{2}{n}S^2-\frac{2r}{n}S.
\endaligned\end{equation}
Comparing $S$ with the corresponding solution of ODE
\begin{equation}\label{nonnegativity2}\aligned
\frac{d}{dt}x=\frac{2}{n}x^2-\frac{2r}{n}x,\ \ \ \  x(0)=S_{\min}(0)
\endaligned\end{equation} gives
\begin{equation}\label{two-eigenvalue3}\aligned
S(x,t)\geq x(t):=\frac{S_{\min}(0)e^{-\frac{2}{n}\int_0^t
r(s)ds}}{1-\frac{2}{n}S_{\min}(0)\int_0^t\Big(e^{-\frac{2}{n}\int_0^s
r(\tilde{s})d\tilde{s}}\Big)ds},
\endaligned\end{equation} where $x(t)$ is the solution of
\eqref{nonnegativity2}. In particular, when $S_{\min}(0)=0$, then
$S(0)\geq0$ and \eqref{two-eigenvalue3} gives
\begin{equation}\label{two-eigenvalue4}\aligned
S(x,t)\geq 0
\endaligned\end{equation} for all $t$. The desired Theorem
\ref{nonnegativity-thm1} is attained.

\vspace*{2mm} \noindent{\bf Proof of Theorem \ref{spec4thm1}.} From
\eqref{4Eig10}, we have
\begin{equation}\label{adds4Eig1}\aligned
\frac{d}{d
t}\Big(\lambda^{-\Delta}e^{\frac{2}{n}\int_0^tr(s)\,ds}\Big)=&e^{\frac{2}{n}\int_0^tr(s)\,ds}\Bigg\{
\int\limits_M\Big(\lambda^{-\Delta}Su^2-S|\nabla u|^2
+2S^{ij}u_{i}u_j\Big)\,dv\Bigg\}\\
\geq &e^{\frac{2}{n}\int_0^tr(s)\,ds}\Bigg\{
\lambda^{-\Delta}\int\limits_M
Su^2\,dv+(2\theta-1)\int\limits_MS|\nabla u|^2 \,dv\Bigg\},
\endaligned\end{equation} which shows (1) holds.
On the other hand, applying \eqref{two-eigenvalue3} into
\eqref{adds4Eig1}, we achieve
\begin{equation}\label{adds4Eig2}\aligned
\frac{d}{d
t}\Big(\lambda^{-\Delta}e^{\frac{2}{n}\int_0^tr(s)\,ds}\Big)\geq
&e^{\frac{2}{n}\int_0^tr(s)\,ds}\Big\{ \lambda^{-\Delta}x(t)
+(2\theta-1)x(t)\lambda^{-\Delta}\Big\}\\
=&2\theta e^{\frac{2}{n}\int_0^tr(s)\,ds}x(t)\lambda^{-\Delta},
\endaligned\end{equation} which shows (2) holds.

\vspace*{2mm} \noindent{\bf Proof of Corollary \ref{spec4thm2}.} As
in \cite{Liyi2014} of Li, using the fact $R_{ij}=\frac{R}{2}g_{ij}$,
we can compute
\begin{equation}\label{adds4Eig9}\aligned
S_{ij}V^iV^j=&\Big(\frac{R}{2}g_{ij}-\alpha\varphi_i\varphi_j\Big)V^iV^j\\
=&\frac{R}{2}|V|^2-\alpha\langle \nabla\varphi,V \rangle^2\\
\geq&\frac{R}{2}|V|^2-\alpha|\nabla\varphi|^2|V|^2\\
\geq&\Big(\frac{R}{2}-\alpha|\nabla\varphi|^2\Big)|V|^2,
\endaligned\end{equation} where $V=(V^i)$. Since $R_{ij}\leq\epsilon
u_iu_j$ and $\epsilon\leq\frac{2\alpha(\theta-1)}{2\theta-1}$ with
$\theta>\frac{1}{2}$, we have
$$\Big(\frac{1}{2}-\theta\Big)R+(\theta-1)\alpha|\nabla\varphi|^2\geq0$$
which is equivalent to
$$\Big(\frac{R}{2}-\alpha|\nabla\varphi|^2\Big)|V|^2\geq\theta
S|V|^2.$$ Therefore, we have $S_{ij}\geq\theta Sg_{ij}$ from
\eqref{adds4Eig9} and the consequence (i) follows.

On the other hand, we can check that if $|\nabla
\varphi|^2g_{ij}\geq2\varphi_i\varphi_j$, then
\begin{equation}\label{secondadds4Eig9}\aligned
S_{ij}V^iV^j=&\frac{R}{2}|V|^2-\alpha\langle \nabla\varphi,V
\rangle^2\\
\geq&\frac{R}{2}|V|^2-\frac{\alpha}{2}|\nabla\varphi|^2|V|^2\\
=&\frac{1}{2}S|V|^2
\endaligned\end{equation} which shows that
$S_{ij}\geq\frac{1}{2}Sg_{ij}$ and the consequence (ii) follows. We
complete the proof of Corollary \ref{spec4thm2}.

\section{Proof of Theorem \ref{twodim4thm3}-Theorem \ref{two4thm5} }

We first prove Theorem \ref{twodim4thm3}.

\vspace*{2mm} \noindent{\bf Proof of Theorem \ref{twodim4thm3}.}
When $n=2$, we have $R_{ij}=\frac{R}{2}g_{ij}$ and
\begin{equation}\label{4Eig14}\aligned
S_{ij}=\frac{R}{2}g_{ij}-\alpha\varphi_i\varphi_j=\frac{1}{2}(S+\alpha|\nabla
\varphi|^2)g_{ij}-\alpha\varphi_i\varphi_j.
\endaligned\end{equation} Hence,
\begin{equation}\label{4Eig15}\aligned
|S_{ij}|^2=&\frac{R^2}{2}+\alpha^2|\nabla \varphi|^4-\alpha R|\nabla
\varphi|^2\\
=&\frac{1}{2}S^2+\frac{1}{2}\alpha^2|\nabla \varphi|^4
\endaligned\end{equation} and \eqref{4Eig9} becomes
\begin{equation}\label{4Eig16}\aligned
\frac{d}{d t}\lambda=&-r\lambda
+\int\limits_M\Big\{2b|S_{ij}|^2u^2+2S^{ij}u_{i}u_j+2b\alpha(\Delta\varphi)^2u^2\\
&+(2b-1)bS^2u^2-(2b-1)\lambda Su^2+(2b-1)S|\nabla u|^2]\Big\}\,dv\\
=&-r\lambda +\int\limits_M\Big\{bS^2u^2+b\alpha^2|\nabla
\varphi|^4u^2+S|\nabla u|^2+\alpha|\nabla \varphi|^2|\nabla u|^2
-2\alpha\langle\nabla
u,\nabla\varphi\rangle^2\\
&+2b\alpha(\Delta\varphi)^2u^2+(2b-1)bS^2u^2-(2b-1)\lambda Su^2+(2b-1)S|\nabla u|^2]\Big\}\,dv\\
=&-r\lambda +\int\limits_M\Big\{2b^2S^2u^2-(2b-1)\lambda
Su^2+2bS|\nabla u|^2\\
&+b\alpha^2|\nabla \varphi|^4u^2+\alpha|\nabla \varphi|^2|\nabla
u|^2 -2\alpha\langle\nabla
u,\nabla\varphi\rangle^2+2b\alpha(\Delta\varphi)^2u^2]\Big\}\,dv.
\endaligned\end{equation}
Therefore, we obtain Theorem \ref{twodim4thm3}

\vspace*{2mm} \noindent{\bf Proof of Theorem \ref{two4thm4}.} In
particular, when $r=0$, the rescaled List's extended Ricci flow
\eqref{Int3} becomes the List's extended Ricci flow \eqref{Int2} and
\eqref{4Eig16} becomes
\begin{equation}\label{two-eigenvalue1}\aligned
\frac{d}{d t}\lambda=&\int\limits_M\Big\{2b^2S^2u^2+(1-2b)\lambda
Su^2+2bS|\nabla u|^2\\
&+b\alpha^2|\nabla \varphi|^4u^2+\alpha|\nabla \varphi|^2|\nabla
u|^2-2\alpha\langle\nabla
u,\nabla\varphi\rangle^2+2b\alpha(\Delta\varphi)^2u^2]\Big\}\,dv,
\endaligned\end{equation} respectively. In particular, for $r=0$, we
have from \eqref{two-eigenvalue3}
\begin{equation}\label{two-eigenvalue4}\aligned
S(x,t)\geq \frac{S_{\min}(0)}{1-tS_{\min}(0)}.
\endaligned\end{equation} Therefore, \eqref{two-eigenvalue1} gives
\begin{equation}\label{two-eigenvalue5}\aligned
\frac{d}{d t}\lambda\geq&\int\limits_M\Big\{2b^2S^2u^2+(1-2b)\lambda
Su^2+2bS|\nabla u|^2\Big\}\,dv\\
\geq&\int\limits_M\Bigg\{2b^2S^2u^2+(1-2b)\lambda
Su^2+\frac{2bS_{\min}(0)}{1-tS_{\min}(0)}|\nabla u|^2\Bigg\}\\
=&\int\limits_M\Bigg\{2b^2S^2u^2+(1-2b)\lambda
Su^2+\frac{2bS_{\min}(0)}{1-tS_{\min}(0)}(\lambda u^2-bSu^2)\Bigg\}\\
\geq&\frac{S_{\min}(0)}{1-tS_{\min}(0)}\lambda
+2b^2\int\limits_M\Big(S^2-S\frac{S_{\min}(0)}{1-tS_{\min}(0)}\Big)u^2\,dv
\endaligned\end{equation}
Using the inequality $y^2-cy\geq-\frac{1}{4}c^2$ into
\eqref{two-eigenvalue5}, we derive
\begin{equation}\label{two-eigenvalue6}\aligned
\frac{d}{d t}\lambda\geq&\frac{S_{\min}(0)}{1-tS_{\min}(0)}\lambda
-\frac{b^2S^2_{\min}(0)}{2[1-tS_{\min}(0)]^2},
\endaligned\end{equation} which shows
\begin{equation}\label{two-eigenvalue7}\aligned
\frac{d}{d
t}\Bigg\{[1-tS_{\min}(0)]\lambda-\frac{b^2S_{\min}(0)}{2}\ln[1-tS_{\min}(0)]\Bigg\}\geq0.
\endaligned\end{equation} Integrating both sides of
\eqref{two-eigenvalue7} on $t$, we have
\begin{equation}\label{two-eigenvalue8}\aligned
\lambda(t)\geq\frac{1}{1-tS_{\min}(0)}\lambda(0)
+\frac{b^2S_{\min}(0)}{2[1-tS_{\min}(0)]}\ln[1-tS_{\min}(0)]
\endaligned\end{equation} and Theorem \ref{two4thm4} follows.

\vspace*{2mm} \noindent{\bf Proof of Theorem \ref{two4thm5}.} Note
that \eqref{4Eig16} becomes
\begin{equation}\label{two-eigenvalue11}\aligned \frac{d}{d
t}\lambda=&-r\lambda +\int\limits_M\Big\{2b^2S^2u^2+(1-2b)\lambda
Su^2+2bS|\nabla u|^2\\
&+b\alpha^2|\nabla \varphi|^4u^2+\alpha|\nabla \varphi|^2|\nabla
u|^2 -2\alpha\langle\nabla
u,\nabla\varphi\rangle^2+2b\alpha(\Delta\varphi)^2u^2]\Big\}\,dv.
\endaligned\end{equation}
In particular, for $n=2$, we have from \eqref{two-eigenvalue3}
\begin{equation}\label{two-eigenvalue14}\aligned
S(x,t)\geq x(t):=\frac{S_{\min}(0)e^{-\int_0^t
r(s)ds}}{1-S_{\min}(0)\int_0^t\Big(e^{-\int_0^s
r(\tilde{s})d\tilde{s}}\Big)ds}.
\endaligned\end{equation}
Hence, \eqref{two-eigenvalue11} yields
\begin{equation}\label{two-eigenvalue15}\aligned \frac{d}{d
t}\lambda\geq&-r\lambda +\int\limits_M\Big\{2b^2S^2u^2+(1-2b)\lambda
Su^2+2bS|\nabla u|^2\Big\}\,dv\\
\geq&-r\lambda +2b^2\int\limits_MS^2u^2\,dv+(1-2b)\lambda
x(t)+2bx(t)\int\limits_M|\nabla u|^2\,dv\\
=&-r\lambda+2b^2\int\limits_MS^2u^2\,dv+(1-2b)\lambda
x(t)+2bx(t)\Big\{\lambda-b\int\limits_MSu^2\,dv\Big\}\\
=&(x(t)-r)\lambda+2b^2\int\limits_M(S-x(t))Su^2\,dv\\
\geq&(x(t)-r)\lambda.
\endaligned\end{equation}
Note that $\lambda(t)>0$. We obtain from \eqref{two-eigenvalue15}
\begin{equation}\label{two-eigenvalue16}\aligned \frac{d}{d
t}(\ln\lambda)\geq x(t)-r.
\endaligned\end{equation}
We complete the proof of Theorem \ref{two4thm5}.

\section{Proof of Theorem \ref{2thm1}-Theorem \ref{3thm1}}

In order to derive our results, we first prove the following three
lemmas:

\begin{lem}\label{2lem1} Under the evolution equation \eqref{2Proof1}, we have
\begin{equation}\label{2Proof4}
(e^{-f}dv)_t=(-f_t-S+r)e^{-f}dv=-(\Delta e^{-f})dv,
\end{equation}
\begin{equation}\label{2Proof6}
(|\nabla f|^2)_t=2S_{ij}f^if^j-\frac{2r}{n}|\nabla f|^2 +2\nabla
f\nabla(-\Delta f+|\nabla f|^2-S+r).
\end{equation}

\end{lem}

\proof The formula \eqref{2Proof4} is a direct conclusion of
\eqref{2eigenvalue2} in Lemma \ref{2eigenvaluelem1} and $f_t=-\Delta
f+|\nabla f|^2-S+r$. Similarly, we have
\begin{equation}\label{add1Lemma5}\aligned
(|\nabla
f|^2)_t=&(g^{ij}f_if_j)_t\\
=&f_if_j\frac{\partial}{\partial t}g^{ij}+2g^{ij}(f_t)_if_j\\
=&2S^{ij}f_if_j-\frac{2r}{n}|\nabla f|^2 +2\nabla f\nabla(-\Delta
f+|\nabla f|^2-S+r)
\endaligned\end{equation} and \eqref{2Proof6} is achieved. We
complete the proof of Lemma \ref{2lem1}.

\begin{lem}\label{2lem2} With the help of Lemma \ref{2lem1}, we have
\begin{equation}\label{2Proof7}
\Big(\int\limits_MS\,e^{-f}dv\Big)_t=\int\limits_M\Big(2|S_{ij}|^2-\frac{2r}{n}S+2\alpha(\Delta
\varphi)^2\Big)\,e^{-f}dv,
\end{equation}
\begin{equation}\label{2Proof8}
\Big(\int\limits_M|\nabla
f|^2\,e^{-f}dv\Big)_t=\int\limits_M\Big(2S_{ij}f^if^j-\frac{2r}{n}|\nabla
f|^2-2\Delta^2f-2\Delta S+\Delta|\nabla f|^2\Big)\,e^{-f}dv,
\end{equation}
\begin{equation}\label{2Proof9}\aligned
\Big(\int\limits_M(|\nabla
f|^2+S)\,e^{-f}dv\Big)_t=&\int\limits_M\Big(2|S_{ij}|^2-\frac{2r}{n}(S+|\nabla
f|^2)+2S_{ij}f^if^j\\
&-2\Delta^2f-2\Delta S+\Delta|\nabla f|^2+2\alpha(\Delta
\varphi)^2\Big)\,e^{-f}dv.
\endaligned\end{equation}

\end{lem}

\proof

By the Stoke formula, we can obtain \eqref{2Proof7} from
\eqref{2eigenvalue3} in Lemma \ref{2eigenvaluelem1}. It is easy to
see that \eqref{2Proof8} holds from \eqref{2Proof6}. Thus, the
formula \eqref{2Proof9} follows by combining \eqref{2Proof7} with
\eqref{2Proof8}.

\begin{lem}\label{2lem3} Under the rescaled List's extended Ricci flow
\eqref{Int3}, we have
\begin{equation}\label{2Proof10}\aligned
\int\limits_M|f_{ij}|^2\,e^{-f}dv=\int\limits_M\Big(\frac{1}{2}\Delta|\nabla
f|^2- \Delta^2 f-S_{ij}f^if^j-\alpha\langle\nabla f,\nabla
\varphi\rangle^2\Big)\,e^{-f}dv,
\endaligned\end{equation}
\begin{equation}\label{2Proof11}\aligned
2\int\limits_MS^{ij}f_{ij}\,e^{-f}dv=\int\limits_M[2S_{ij}f^if^j-\Delta
S+2\alpha(\Delta\varphi)\langle\nabla
f,\nabla\varphi\rangle]\,e^{-f}dv,
\endaligned\end{equation}
\begin{equation}\label{2Proof12}\aligned
\int\limits_M|S_{ij}+f_{ij}|^2\,e^{-f}dv=&\int\limits_M[|S_{ij}|^2+S_{ij}f^if^j+\frac{1}{2}\Delta|\nabla
f|^2- \Delta^2 f-\Delta S\\
&-\alpha\langle\nabla f,\nabla
\varphi\rangle^2+2\alpha(\Delta\varphi)\langle\nabla
f,\nabla\varphi\rangle]\,e^{-f}dv.
\endaligned\end{equation}

\end{lem}

\proof  Multiplying both sides of the following well-known Bochner
formula with $e^{-f}$
\begin{equation}\label{2Proof13}
\frac{1}{2}\Delta|\nabla f|^2=|f_{ij}|^2+2\nabla f\nabla\Delta
f+R_{ij}f^if^j\end{equation} and integrating on it, we have
\begin{equation}\label{2Proof14}\aligned
\int\limits_M|f_{ij}|^2\,e^{-f}dv=&\int\limits_M\Big(\frac{1}{2}\Delta|\nabla
f|^2- \Delta^2 f-R_{ij}f^if^j\Big)\,e^{-f}dv\\
=&\int\limits_M\Big(\frac{1}{2}\Delta|\nabla f|^2- \Delta^2
f-S_{ij}f^if^j-\alpha\langle\nabla f,\nabla
\varphi\rangle^2\Big)\,e^{-f}dv
\endaligned\end{equation} and \eqref{2Proof10} follows.
On the other hand,
\begin{equation}\label{2Proof15}\aligned
2\int\limits_MS_{ij}f^if^j\,e^{-f}dv=&2\int\limits_M(S^{ij}f_i)_{,j}\,e^{-f}dv\\
=&2\int\limits_M(S^{ij}{}_{,j}f_i+S^{ij}f_{ij})\,e^{-f}dv.
\endaligned\end{equation}
Applying \eqref{2Proof10} into \eqref{2Proof15} gives
\begin{equation}\label{2Proof16}\aligned
2\int\limits_MS^{ij}f_if_j\,e^{-f}dv
=&2\int\limits_M\Big(\frac{1}{2}S_{,i}f^i-\alpha(\Delta\varphi)\varphi_if^i+S^{ij}f_{ij}\Big)\,e^{-f}dv\\
=&\int\limits_M[\Delta S-2\alpha(\Delta\varphi)\langle\nabla
f,\nabla\varphi\rangle+2S^{ij}f_{ij}]\,e^{-f}dv,
\endaligned\end{equation} which gives the desired formula
\eqref{2Proof11}.

By combining \eqref{2Proof10} with \eqref{2Proof11}, we derive
\eqref{2Proof12} finally.


\vspace*{2mm}
\noindent{\bf Proof of Theorem \ref{2thm1}.}
Applying\eqref{2Proof12} into \eqref{2Proof9}, we obtain
\begin{equation}\label{2Proof19}\aligned
\frac{d}{dt}\mathcal{F}=&\frac{d}{dt}\int\limits_M(|\nabla
f|^2+S)\,e^{-f}dv\\
=&\int\limits_M\Big(2|S_{ij}|^2-\frac{2r}{n}(S+|\nabla
f|^2)+2S^{ij}f_if_j\\
&-2\Delta^2f-2\Delta S+\Delta|\nabla f|^2+2\alpha(\Delta
\varphi)^2\Big)\,e^{-f}dv\\
=&\int\limits_M\Big(-\frac{2r}{n}(|\nabla
f|^2+S)+2|S_{ij}+f_{ij}|^2\\
&+2\alpha(\Delta \varphi)^2+2\alpha\langle\nabla f,\nabla
\varphi\rangle^2-4\alpha(\Delta\varphi)\langle\nabla
f,\nabla\varphi\rangle\Big)\,e^{-f}dv\\
=&-\frac{2r}{n}\mathcal{F}+2\int\limits_M\Big(|S_{ij}+f_{ij}|^2+\alpha|\Delta
\varphi-\langle\nabla f,\nabla\varphi\rangle|^2\Big)\,e^{-f}dv.
\endaligned\end{equation}
Therefore, we obtain
\begin{equation}\label{2Proof20}\aligned
\frac{d}{dt}\mathcal{F}_k=&\frac{d}{dt}\mathcal{F}+(k-1)\frac{d}{dt}
\int\limits_MS\,e^{-f}dv\\
=&-\frac{2r}{n}\mathcal{F}_k+2(k-1)\int\limits_M\Big(|S_{ij}|^2+\alpha(\Delta
\varphi)^2\Big)\,e^{-f}dv\\
&+2\int\limits_M\Big(|S_{ij}+f_{ij}|^2+\alpha|\Delta
\varphi-\langle\nabla f,\nabla\varphi\rangle|^2\Big)\,e^{-f}dv
\endaligned\end{equation} and the desired formula \eqref{add2Proof17} is achieved.

The formula \eqref{add2Proof18} can be achieved by a direct
computation.

\vspace*{2mm} \noindent{\bf Proof of Theorem \ref{1th1}.} In order
to prove Theorem \ref{1th1}, we fist give the following key
proposition:

\begin{prop}\label{2thm2}

Let $(M^n,g(t))$ be a compact Riemannian manifold with $g(t)$
satisfying the rescaled List's extended Ricci flow \eqref{Int3}. We
let $\lambda(t)$ be the lowest eigenvalue of the operator
$-4\Delta+kS$ with $k\geq1$. If there exists a function $v=v(x,t)$
such that
\begin{equation}\label{2Proof21}\aligned
r(t)=\frac{\int_M(|\nabla v|^2+kS)\,e^{-v}dv}{k\int_M\,e^{-v}dv},
\endaligned\end{equation}
then the lowest eigenvalue eigenvalue $\lambda(t)$ is nondecreasing
under \eqref{Int3} provided $r(t)\leq0$. The monotonicity is strict
unless the metric is Einstein.

\end{prop}

\proof Since the lowest eigenvalue of the operator $-4\Delta+kS$ on
a compact Riemannian manifold is given by
\begin{equation}\label{2Proof22}\aligned
\lambda(g(t))=\inf\limits_f
\int\limits_M\mathcal{F}_k(g,f)\,e^{-f}dv,
\endaligned\end{equation} where the infimum is taken over functions
satisfying $\int_M\,e^{-f}dv=1$. For the compact Riemannian
manifold, the lowest eigenvalue $\lambda$ can be attained by a
smooth function $u$. Therefore, there exists a smooth function $u$
such that $(-4\Delta+kS)u=\lambda u$ and
\begin{equation}\label{2Proof22}\aligned
\lambda(g(t))=\int\limits_M(4|\nabla
u|^2+kSu^2)dv=\int\limits_M(|\nabla
\tilde{f}|^2+kS)\,e^{-\tilde{f}}dv
\endaligned\end{equation}
by letting $-\tilde{f}=2\ln u$. Applying \eqref{2Proof22} into
\eqref{add2Proof18} gives
\begin{equation}\label{2Proof23}\aligned
\frac{d}{dt}\lambda=&\frac{2r}{n}(\lambda-kr)
+2(k-1)\int\limits_M\Big(|S_{ij}-\frac{r}{n}g_{ij}|^2+\alpha(\Delta
\varphi)^2\Big)\,e^{-\tilde{f}}dv\\
&+2\int\limits_M\Big(|S_{ij}+\tilde{f}_{ij}-\frac{r}{n}g_{ij}|^2+\alpha|\Delta
\varphi-\langle\nabla
\tilde{f},\nabla\varphi\rangle|^2\Big)\,e^{-\tilde{f}}dv.
\endaligned\end{equation}
By using the assumption \eqref{2Proof21}, we get
\begin{equation}\label{2Proof24}\aligned
kr=&\int\limits_M(|\nabla v|^2+kS)\,e^{-v}dv\\
\geq&\int\limits_M(|\nabla
\tilde{f}|^2+kS)\,e^{-\tilde{f}}dv\\
=&\lambda
\endaligned\end{equation} which means that $\lambda-kr\leq0$.
Therefore, if $r(t)\leq0$, then \eqref{2Proof23} shows that
$\lambda(t)$ is nondecreasing under \eqref{Int3}. The monotonicity
is strict unless
\begin{equation}\label{2Proof25}\aligned
S_{ij}-\frac{r}{n}g_{ij}=0,\ \  \ \Delta \varphi=0
\endaligned\end{equation} and
\begin{equation}\label{2Proof26}\aligned
S_{ij}+\tilde{f}_{ij}-\frac{r}{n}g_{ij}=0,\ \ \ \Delta
\varphi-\langle\nabla \tilde{f},\nabla\varphi\rangle=0.
\endaligned\end{equation}
Notice that the Riemannian manifold is compact. Hence, we have
$\varphi$ is constant from the second equality in \eqref{2Proof25}
or the second equality in \eqref{2Proof26}. Therefore, the metric
$g$ is Einstein. We complete the proof of Proposition \ref{2thm2}.

Now, we are in a position to prove Theorem\ref{1th1}.  We note that
for the extended Hamilton normalized flow under the List's extended
Ricci flow, we have $$r(t)=\frac{\int_MS\,dv}{\int_M\,dv}.$$
Choosing $v=\ln( {\rm Vol}(M^n))$ in \eqref{2Proof21}, we derive
Theorem \ref{1th1}.


\vspace*{2mm} \noindent{\bf Proof of Theorem \ref{3thm1}.} Since the
$\mathcal{W}$-functional is related with $\mathcal{F}$ by
\begin{equation}\label{3Ent3}\aligned
\mathcal{W}=&\tau^2\int\limits_M(S+\frac{n}{2\tau}+|\nabla
f|^2)e^{-f}\,dv\\
=&\tau^2\mathcal{F}+\frac{n\tau}{2},
\endaligned\end{equation} with the help of \eqref{add2Proof17}, we
obtain
\begin{equation}\label{3Ent4}\aligned
\frac{d}{d t}\mathcal{W}=&\tau^2\frac{d}{d
t}\mathcal{F}+2\tau\mathcal{F}+\frac{n}{2}\\
=&\tau^2\Bigg\{-\frac{2r}{n}\mathcal{F}+2\int\limits_M\Big(|S_{ij}+f_{ij}|^2+\alpha|\Delta
\varphi-\langle\nabla
f,\nabla\varphi\rangle|^2\Big)\,e^{-f}dv\Bigg\}\\
&+2\tau\mathcal{F}+\frac{n}{2}.
\endaligned\end{equation}
Applying
$$
2\tau^2\int\limits_M|S_{ij}+f_{ij}+\frac{1}{2\tau}g_{ij}|^2\,e^{-f}dv
=2\tau^2\int\limits_M|S_{ij}+f_{ij}|^2\,e^{-f}dv+2\tau\mathcal{F}+\frac{n}{2}
$$
into \eqref{3Ent4} yields
\begin{equation}\label{3Ent6}\aligned
\frac{d}{d t}\mathcal{W}=&2\tau^2\Bigg\{-\frac{r}{n}\mathcal{F}
+\int\limits_M\Big(|S_{ij}+f_{ij}+\frac{1}{2\tau}g_{ij}|^2\\
&+\alpha|\Delta \varphi-\langle\nabla
f,\nabla\varphi\rangle|^2\Big)\,e^{-f}dv\Bigg\}.
\endaligned\end{equation}
By the definition of $\mathcal{W}_k$, we know
\begin{equation}\label{3Ent7}\aligned
\mathcal{W}_k=\mathcal{W}+(k-1)\tau^2\int\limits_MSe^{-f}\,dv+\frac{(k-1)n\tau}{2}.
\endaligned\end{equation}
Thus, we get from \eqref{2Proof7} and \eqref{3Ent6}
\begin{equation}\label{3Ent8}\aligned
\frac{d}{d t}\mathcal{W}_k=&\frac{d}{d
t}\mathcal{W}+\frac{(k-1)n}{2}+2(k-1)\tau\int\limits_MSe^{-f}\,dv\\
&+(k-1)\tau^2\Bigg(\int\limits_MSe^{-f}\,dv\Bigg)_t\\
=&2\tau^2\Bigg\{-\frac{r}{n}\mathcal{F}
+\int\limits_M\Big(|S_{ij}+f_{ij}+\frac{1}{2\tau}g_{ij}|^2\\
&+\alpha|\Delta \varphi-\langle\nabla
f,\nabla\varphi\rangle|^2\Big)\,e^{-f}dv\Bigg\}+\frac{(k-1)n}{2}+2(k-1)\tau\int\limits_MSe^{-f}\,dv\\
&+(k-1)\tau^2\int\limits_M\Big(2|S_{ij}|^2-\frac{2r}{n}S+2\alpha(\Delta
\varphi)^2\Big)\,e^{-f}dv\\
=&2\tau^2\Bigg\{-\frac{r}{n}\mathcal{F}_k+(k-1)\int\limits_M\Big(|S_{ij}+\frac{1}{2\tau}g_{ij}|^2
+\alpha(\Delta \varphi)^2\Big)\,e^{-f}dv\\
&+\int\limits_M\Big(|S_{ij}+f_{ij}+\frac{1}{2\tau}g_{ij}|^2
+\alpha|\Delta \varphi-\langle\nabla
f,\nabla\varphi\rangle|^2\Big)\,e^{-f}dv\Bigg\},
\endaligned\end{equation} and the desired formula \eqref{3Ent2}
is attained.

The formula \eqref{add3Ent2} can be achieved by a direct
computation.

\bibliographystyle{Plain}

\end{document}